\documentclass[a4paper]{article}
\usepackage{RR}
\usepackage{hyperref}
\usepackage[utf8]{inputenc}
\usepackage[english]{babel}

\RRdate{Juin 2010}

\RRauthor{
Robin Genuer 
  \thanks[sfn]{Univ Paris-Sud, Laboratoire de Math\'ematique, UMR 8628, Orsay F-91405}
\thanks{Inria Saclay Ile-de-France}
}
\authorhead{Genuer}
\RRtitle{Bornes de risque pour les forêts purement uniformément aléatoires.}
\RRetitle{Risk bounds for purely uniformly random forests}
\titlehead{Risk bounds for purely uniformly random forests}
\RRresume{Introduites par Leo Breiman en 2001, les forêts aléatoires sont une méthode statistique très performante. D'un point de vue théorique, leur analyse est difficile, du fait de la complexité de l'algorithme. Pour expliquer ces performances, des versions de forêts aléatoires simplifiées, et donc plus faciles à analyser, ont été introduites. Ces versions ont été appelées forêts purement aléatoires. Dans cet article, nous introduisons une autre version simplifiée, que nous appelons forêts purement uniformément aléatoires. Dans un contexte de régression, avec une seule variable explicative, nous montrons que les arbres aléatoires ainsi que les forêts aléatoires atteignent la vitesse de convergence minimax. De plus, nous prouvons que les forêts aléatoires améliorent les performances des arbres aléatoires, en réduisant la variance des estimateurs associés d'un facteur de trois quarts.
}
\RRabstract{Random forests, introduced by Leo Breiman in 2001, are a very effective statistical method. The complex mechanism of the method makes theoretical analysis difficult. Therefore, a simplified version of random forests, called purely random forests, which can be theoretically handled more easily, has been considered. In this paper we introduce a variant of this kind of random forests, that we call purely uniformly random forests. In the context of regression problems with a one-dimensional predictor space, we show that both random trees and random forests reach minimax rate of convergence. In addition, we prove that compared to random trees, random forests improve accuracy by reducing the estimator variance by a factor of three fourths.
}
\RRmotcle{{\sc For\^ets al\'eatoires, R\'egression non-param\'etrique, Vitesse de convergence, Randomisation.}}
\RRkeyword{{\sc Random Forests, Non-parametric regression, Rate of convergence, Randomization.}}
\RRprojets{{\sc Select}}
\RRtheme{\THCog} 
 \RCSaclay 

\usepackage{amssymb,amsmath}
\usepackage{dsfont}
\usepackage{amsmath}

\everymath{\displaystyle}

\newcommand{\E}{\mathbb{E}}
\newcommand{\C}{\mathbb{C}}

\renewcommand{\P}{\textup{P}}

\newcommand{\var}{\textup{Var}}
\newcommand{\eps}{\varepsilon}

\renewcommand{\leq}{\leqslant}

 \newcommand{\petito}[1]{\underset{#1}{\mathrm{o}}}

\newtheorem{prop}{Proposition}
\newtheorem{cor}{Corollary}
\newtheorem{theo}{Theorem}

\newtheorem{lem}{Lemma}
\newtheorem{rem}{Remark}

\begin{document}
\RRNo{7318}
\makeRR   

\section{Introduction}

Random forests (RF), introduced by Leo Breiman in 2001 \cite{Breiman01}, are a very effective statistical method. They give outstanding performances in a lot of situations for both regression and classification problems. Mathematical understanding of these good performances remains quite unknown. As defined by Leo Breiman, a random forest is a collection of tree-predictors $\{h(x,\Theta_l),1\leq l \leq q\}$, where $(\Theta_l)_{1\leq l \leq q}$ are i.i.d. random vectors, and a random forest predictor is obtained by aggregating this collection of trees. In addition to consistency results, one of the main theoretical challenges is to explain why a random forest improves so much the performance of a single tree.

In \cite{Breiman01}, Leo Breiman introduced a specific instance of random forest, called random forests-RI, which has been adopted in many fields as a reference method. Indeed, random forests-RI are simple to use, and are efficiently coded in the popular R-package \texttt{randomForest} \cite{Liaw02}. They are effective for a predictive goal and they can also be used for variable selection (see e.g. \cite{Diaz06}, \cite{Genuer10}).

However, forests-RI are very difficult to handle theoretically. This is why people are interested in simplified versions, called purely random forests (PRF). The main difference is that in PRF, the splits of tree nodes are randomly drawn \textit{independently} of the learning sample; while in random forests-RI, the splits are optimized using the learning sample. This independence between splits and learning sample makes mathematical analysis easier. In \cite{Cutler01}, Cutler and Zhao introduced PERT (Perfect Random Tree Ensemble), an algorithm which builds some purely random forests, and illustrated its good performance on benchmark datasets. More recently Biau et al. \cite{Biau08} showed that both purely random trees and purely random forests are universally consistent.

Our paper offers to examine another simple variant of random forests, which can be put in the so-called purely random forests family. We call it \textit{purely uniformly random forests} and we analyze its risk, only in a regression framework with a one-dimensional predictor space. The main goal is to emphasize the gain of using a forest instead of a tree. The results of this paper are twofold: first we show that both purely uniformly random trees and forests risks reach minimax rate of convergence on the Lipschitz functions class; second we show that forests improve the variance term by a factor of three fourths while not increasing the bias.

The paper is organized as follows. Section~\ref{framework} presents the model. Section~\ref{risk-tree} and Section~\ref{risk-forest} give some risk bounds for purely uniformly random trees and purely uniformly random forests respectively. Section~\ref{conclusion} concludes the paper, while proofs are collected in Section~\ref{proofs}.

\section{Framework}\label{framework}
The framework we consider all along the paper is the classical random design regression framework.

More precisely, consider a learning set $\mathcal{L}_n=\{(X_1,Y_1), \ldots, (X_n,Y_n)\}$ made of $n$ i.i.d. observations of a vector $(X,Y)$ from an unknown distribution. $Y$ is real-valued since we are in a regression framework. We consider the following statistical model:
\begin{align}\label{model}
 Y_i=s(X_i)+\eps_i \quad \mbox{for }i=1,\ldots,n \; .
\end{align}
$s$ is the unknown regression function and the goal is to estimate $s$. We make the following assumptions on model (\ref{model}):
\textit{
\begin{itemize}
 \item $X\in [0,1]$ with continuous density function $ \mu $;
 \item $(\eps_1,\ldots,\eps_n)$ are i.i.d. observations of $\eps$, independent of $\mathcal{L}_n$, with $\E[\eps]=0$ and where $\var(\eps)=\sigma^2$ is assumed to be known.
\end{itemize}
}
Note that we deal only with a one-dimensional predictor space.

This paper aims at comparing performances in estimating $s$ using a single random tree and a random forest of a special kind, described in the next section.

\section{Risk bounds for Purely Uniformly Random Trees}\label{risk-tree}

\subsection{Tree definition}\label{def-tree}

The principle of Purely Uniformly Random Trees (PURT) is that we draw $k$ uniform random variables, which form the partition of the input space $[0,1]$. Then we build a regressogram on this partition, that we call a tree.

Note that, unlike purely random forests or random forests-RI, the tree structure of individual predictors is not obvious. This comes from the fact that in PURT the partition is not obtained in a recursive manner. Nevertheless we keep the vocabulary of trees and forests to distinguish individual predictors from aggregated ones.

Let us mention that, all along the paper, we make a slight language abuse. Indeed, we refer to random tree, the tree himself (as a graph), the corresponding partition of $[0,1]$, as well as the corresponding estimator.

\bigskip

More precisely, let $\mathbb{U}=(U_1,\ldots,U_k)$ be $k$ i.i.d. random variables of uniform distribution on $[0,1]$, where $k$ is a natural integer which will depend on the number of observations $n$.

A Purely Uniformly Random Tree (PURT), associated with $\mathbb{U}$, is defined for $x\in[0,1]$ as:
$$\hat{s}_{\mathbb{U}}(x) = \sum_{j=0}^k \hat{\beta}_j \mathds{1}_{U_{(j)} < x \leq U_{(j+1)}}$$
where
$$\hat{\beta}_j=\frac{1}{\sharp\{i: \; U_{(j)} < X_i \leq U_{(j+1)}\}} \sum_{i: \; U_{(j)} < X_i \leq U_{(j+1)}} Y_i$$
and $(U_{(1)},\ldots,U_{(k)})$ is the ordered statistics of $(U_1,\ldots,U_k)$ and $U_{(0)}=0$, $U_{(k+1)}=1$. $\sharp \, \mathcal{E}$ denotes the cardinality of the set $ \mathcal{E} $.

\begin{rem}
Let us mention that if $\: \sharp\{i: \; U_{(j)} < X_i \leq U_{(j+1)}\} = 0 $, we set $\hat{\beta}_j = 0 $. However as we will see in Section~\ref{tree-var}, our assumptions on $k$ and $n$ will make the probability of observing such an event tend to $0$.
\end{rem}

In addition, let us define, for $x\in[0,1]$:

$$\tilde{s}_{\mathbb{U}}(x) = \sum_{j=0}^k \beta_j \mathds{1}_{U_{(j)} < x \leq U_{(j+1)}}$$
where
$$\beta_j=\E[Y \, | \, U_{(j)} < X \leq U_{(j+1)}] \; .$$
Conditionally on $\mathbb{U}$, $\tilde{s}_{\mathbb{U}}$ is the best approximation of $s$ among all the regressograms based on $\mathbb{U}$, but of course it depends on the unknown distribution of $(X,Y)$.

With these notations, we can write a bias-variance decomposition of the quadratic risk of $\hat{s}_{\mathbb{U}}$ as follows:

\begin{align}\label{decomposition}
\E[(\hat{s}_{\mathbb{U}}(X) - s(X))^2] &= \E[(\hat{s}_{\mathbb{U}}(X) - \tilde{s}_{\mathbb{U}}(X))^2] + \E[(\tilde{s}_{\mathbb{U}}(X) - s(X))^2] \\
&= \mbox{variance term} + \mbox{bias term} \notag
\end{align}

\medskip

\noindent To clarify these variance and bias terms, we emphasize that for a given partition $u$ and a given $x$, we have
$$ \E[\hat{s}_{u}(x)] = \tilde{s}_{u}(x) $$
so $\E[(\hat{s}_{u}(x) - \tilde{s}_{u}(x))^2]$ is the variance of the estimator $\hat{s}_{u}(x)$ and $\E[(\tilde{s}_{u}(x) - s(x))^2]$ is its bias. We then integrate with respect to (w.r.t) $X$ and $ \mathbb{U} $ to get decomposition (\ref{decomposition}).

\subsection{Variance of a tree}\label{tree-var}

We start to deal with the variance term of decomposition~(\ref{decomposition}). First, we work conditionally on $\mathbb{U}$, then the problem reduces to the case of a regressogram on a deterministic partition, and we can apply the following proposition which comes from Arlot \cite{Arlot}.
\begin{prop}\label{Arlot}
Conditionally on $\mathbb{U}$, the variance term of decomposition~(\ref{decomposition}) satisfies:
\begin{equation}\label{var-cond}
\E[(\hat{s}_{\mathbb{U}}(X) - \tilde{s}_{\mathbb{U}}(X))^2 \, | \, \mathbb{U}] = \frac{1}{n} \sum_{j=0}^k (1+\delta_{n,p_j})(\sigma^2 + (\sigma_j^d)^2)
\end{equation}
where
\begin{itemize}
 \item $p_j=\P(U_{(j)} < X \leq U_{(j+1)})$,
 \item $(\sigma_j^d)^2=\E[(s(X) - \tilde{s}_{\mathbb{U}}(X))^2 \, | \, U_{(j)} < X \leq U_{(j+1)}]$,
 \item $\delta_{n,p} \xrightarrow[np \to +\infty]{} 0$.
\end{itemize}
\begin{flushright}$ \blacksquare $\end{flushright}
\end{prop}
We now integrate equation~(\ref{var-cond}) w.r.t. $\mathbb{U}$, and we get the following equality:
\begin{align}\label{var-tree-gen}
 \E[(\hat{s}_{\mathbb{U}}(X) - \tilde{s}_{\mathbb{U}}(X))^2] = \frac{1}{n}
\sum_{j=0}^k \left( \sigma^2 + \sigma^2\E[\delta_{n,p_j}] + \E[(\sigma_j^d)^2] + \E[(\sigma_j^d)^2\delta_{n,p_j}]\right)
\end{align}
Let us stress that equation~(\ref{var-tree-gen}) is general, since it does not depend on the distribution of $\mathbb{U}$. Hence, it can be used for any random partition distributions.

Finally, using the fact that, in our case, $\mathbb{U}$ is made of $k$ i.i.d. random variables of uniform distribution on $[0,1]$, we deduce from equation~(\ref{var-tree-gen}) the following proposition:
\begin{prop}\label{prop-var-tree}
If $k \xrightarrow[n \to +\infty]{} +\infty$, $\frac{k}{n} \xrightarrow[n \to +\infty]{} 0$, $ \mu>0$ and $s$ is $C$-Lipschitz, the variance  of a PUR Tree satisfies:
\begin{equation}\label{var-purt}
 \E[(\hat{s}_{\mathbb{U}}(X) - \tilde{s}_{\mathbb{U}}(X))^2] = \frac{\sigma^2 (k+1)}{n} + \petito{n\to+\infty} \left(\frac{k}{n} \right)
\end{equation}
where the notation $\petito{n\to+\infty} \left(\frac{k}{n} \right)$ denotes a function $f(n)$ such as $ \: \frac{f(n)}{k/n} \xrightarrow[n \to +\infty]{} 0 $.
\begin{flushright}$ \blacksquare $\end{flushright}
\end{prop}
Details of the proof of Proposition~\ref{prop-var-tree} can be found in Section~\ref{proof-prop-var-tree}.

The first two hypotheses of Proposition~\ref{prop-var-tree} ($k \xrightarrow[n \to +\infty]{} +\infty$, $\frac{k}{n} \xrightarrow[n \to +\infty]{} 0$) are the same natural conditions found by Biau et al. \cite{Biau08} for consistency of PRF. They guarantee that the number of splits of the tree must grow to infinity but slower than the number of samples.


\subsection{Bias of a tree}

We now turn to the bias term of decomposition~(\ref{decomposition}). Direct calculations (see Section~\ref{proof-prop-bias-tree} for details) lead to the following upper bound for the bias term of a PURT:
\begin{prop}\label{prop-bias-tree}
If $\mu$ is bounded by $M>0$ and $s$ is $C$-Lipschitz, the bias of a PURT is upper bounded by:
\begin{equation}\label{bias-purt}
 \E[(\tilde{s}_{\mathbb{U}}(X) - s(X))^2] \leq \frac{6MC^2}{(k+1)^2}
\end{equation}
\begin{flushright}$ \blacksquare $\end{flushright}
\end{prop}

\subsection{Risk bounds for a tree}

Putting together (\ref{var-purt}) and (\ref{bias-purt}) leads to the following risk bound for a PURT.
\begin{theo}\label{theo-risk-bound-tree}
If $k \xrightarrow[n \to +\infty]{} +\infty$, $\frac{k}{n} \xrightarrow[n \to +\infty]{} 0$, $0 < \mu \leq M$ and $s$ is $C$-Lipschitz, the risk of a PURT satisfies:
\begin{equation}\label{risk-bound-purt}
 \E[(\hat{s}_{\mathbb{U}}(X) - s(X))^2] \leq \frac{\sigma^2 (k+1)}{n} + \frac{6MC^2}{(k+1)^2} + \petito{n \to +\infty} \left( \frac{k}{n} \right)
\end{equation}
\begin{flushright}$ \blacksquare $\end{flushright}
\end{theo}

The balance between the two first terms of the right hand side (r.h.s.) of (\ref{risk-bound-purt}) leads to take $(k+1)=n^{1/3}$, and gives the following upper bound for the risk of a PURT.
\begin{cor} Under the assumptions of Theorem~\ref{theo-risk-bound-tree},
\begin{equation*}
 \E[(\hat{s}_{\mathbb{U}}(X) - s(X))^2]  \leq K n^{-2/3}  + \petito{n \to +\infty} (n^{-2/3})
\end{equation*}
where $K$ is a positive constant.
\begin{flushright}$ \blacksquare $\end{flushright}
\end{cor}
Therefore, a PURT reaches the minimax rate of convergence associated with the class of Lipschitz functions (see e.g. Ibragimov and Khasminskii \cite{Ibragimov}).

\medskip

Let us now analyze purely uniformly random forests. As a result, we emphasize an improvement given by a forest compared to a single tree.

\section{Risk bounds for Purely Uniformly Random Forests}\label{risk-forest}

\subsection{Forest definition}

A random forest is the aggregation of a collection of random trees. So, in the context of Purely Uniformly Random Forests (PURF), the principle is to generate several PUR Trees by drawing several random partitions given by uniform random variables, and to aggregate them.

\bigskip

Let $\mathbb{V}=(\mathbb{U}^1,\ldots,\mathbb{U}^q)$ be $q$ i.i.d. random vectors of the same distribution as $\mathbb{U}$ (defined in Section~\ref{def-tree}). That is for $l=1,\ldots,q, \quad \mathbb{U}^l = (U_1^l,\ldots,U_k^l)$ where the $ (U_j^l)_{1\leq j \leq k} $ are i.i.d. random variables of uniform distribution on $[0,1]$.

A PURF, associated with $\mathbb{V}$, is defined for $x\in[0,1]$ as follows:
$$\hat{s}(x) = \frac{1}{q} \sum_{l=1}^q \hat{s}_{\mathbb{U}^l}(x) \; .$$
Let us define, for $x\in[0,1]$:
$$\tilde{s}(x) = \frac{1}{q} \sum_{l=1}^q \tilde{s}_{\mathbb{U}^l}(x) \; .$$
Again, we have a bias-variance decomposition of the quadratic risk of $\hat{s}$, given by:
\begin{align}\label{decomposition-forest}
\E[(\hat{s}(X) - s(X))^2] &= \E[(\hat{s}(X) - \tilde{s}(X))^2] + \E[(\tilde{s}(X) - s(X))^2]  \\
 &= \mbox{variance term} + \mbox{bias term} \notag
\end{align}

\subsection{Variance of a forest}

We first deal with the variance term of decomposition~(\ref{decomposition-forest}). We begin to show that when letting the number of trees $q$ grow to infinity, the variance of a PURF is close to the covariance between two PURT.

Indeed, since $\hat{s}(x) = \frac{1}{q} \sum_{l=1}^q \hat{s}_{\mathbb{U}^l}(x)$, the variance term satisfies:

\begin{multline*}
\E[(\hat{s}(X) - \tilde{s}(X))^2] = \frac{1}{q^2} \sum_{l=1}^q \E[(\hat{s}_{\mathbb{U}^l}(X) - \tilde{s}_{\mathbb{U}^l}(X))^2] \\ + \frac{1}{q^2} \sum_{l \neq q} \E[(\hat{s}_{\mathbb{U}^l}(X) - \tilde{s}_{\mathbb{U}^l}(X))(\hat{s}_{\mathbb{U}^m}(X) - \tilde{s}_{\mathbb{U}^m}(X))]
\end{multline*}
\begin{multline*}
 = \frac{1}{q} \E[(\hat{s}_{\mathbb{U}^1}(X) - \tilde{s}_{\mathbb{U}^1}(X))^2] \\ + \frac{q(q-1)}{q^2} \E[(\hat{s}_{\mathbb{U}^1}(X) - \tilde{s}_{\mathbb{U}^1}(X))(\hat{s}_{\mathbb{U}^2}(X) - \tilde{s}_{\mathbb{U}^2}(X))]
\end{multline*}
%
\noindent where the last equality comes from the fact that the $\left( (\hat{s}_{\mathbb{U}^l}(X) -\tilde{s}_{\mathbb{U}^l}(X)\right)_{1\leq l \leq q} $ are of the same distribution.

Now, if we let $ q $ grow to infinity, we get:
\begin{equation*}
 \E[(\hat{s}(X) - \tilde{s}(X))^2] = \E[(\hat{s}_{\mathbb{U}^1}(X) - \tilde{s}_{\mathbb{U}^1}(X))(\hat{s}_{\mathbb{U}^2}(X) - \tilde{s}_{\mathbb{U}^2}(X))] + \petito{q\to+\infty}(1)
\end{equation*}

The next step is to upper bound the covariance between two PURT
$$\E[(\hat{s}_{\mathbb{U}^1}(X) - \tilde{s}_{\mathbb{U}^1}(X))(\hat{s}_{\mathbb{U}^2}(X) - \tilde{s}_{\mathbb{U}^2}(X))] $$
(it is detailed in Section~\ref{proof-theo-var-forest}) and it leads to the following theorem, which gives the behavior of the variance of a PURF:
\begin{theo}\label{theo-var-forest}
If $k \xrightarrow[n \to +\infty]{} +\infty$, $\frac{k}{n} \xrightarrow[n \to +\infty]{} 0$, $\mu>0$, $s$ is $C$-Lipschitz and $ q \xrightarrow[n \to +\infty]{} +\infty $, the variance of a PURF satisfies the following upper bound:
\begin{equation}\label{var-purf}
 \E[(\hat{s}(X) - \tilde{s}(X))^2] \leq \frac{3}{4}\frac{\sigma^2(k+1)}{n} + \petito{n\to+\infty} \left( \frac{k}{n} \right)
\end{equation}
\begin{flushright}$ \blacksquare $\end{flushright}
\end{theo}
Theorem~\ref{theo-var-forest} is to be compared with Proposition~\ref{prop-var-tree} and tells us that the variance of a PUR Forest is upper bounded by three fourths times the variance of a PUR Tree. So, the rate of decay (in terms of power of $n$) of the PUR Forest variance is the same as the PUR Tree variance, and the actual gain appears in the multiplicative constant.

We mention that, as in the analysis of the variance of a tree (see equation~(\ref{var-tree-gen})), we derive, in the proof of Theorem~\ref{theo-var-forest}, a general statement (see equation~(\ref{prop-var-forest-gen}) in Section~\ref{proof-theo-var-forest}), which does not depend on the distribution of the partition defining the random trees.

Let us, finally, comment the hypotheses of Theorem~\ref{theo-var-forest}. First, note that the hypotheses on $k$ and $n$ are the same as in Proposition~\ref{prop-var-tree}, which allows a fair comparison between the two results.
Second, the hypothesis on $q$ allows to ensure that the upper bound on the covariance (given by Corollary~\ref{coro} in Section~\ref{proof-theo-var-forest}) leads to the same upper bound for the variance of the forest.
Finally, the other hypotheses ($\mu>0$, $s$ is $C$-Lipschitz) are the same as in Proposition~\ref{prop-var-tree} and help to control negligible terms.

\subsection{Bias of a forest}

We now deal with the bias term of decomposition~(\ref{decomposition-forest}). A convex inequality gives that the bias of a forest is not larger than the bias of a single tree:
\begin{align*}
\E[(\tilde{s}(X) - s(X))^2] & \leq  \frac{1}{q} \sum_{l=1}^q \E[(\tilde{s}_{\mathbb{U}^l}(X) - s(X))^2] \\
& =  \E[(\tilde{s}_{\mathbb{U}^1}(X) - s(X))^2] \; .
\end{align*}
So from Proposition~\ref{prop-bias-tree}, we deduce that:
\begin{prop}
If $\mu$ is bounded by $M>0$ and $s$ is $C$-Lipschitz, the bias of a PURF satisfies the same inequality as (\ref{bias-purt}), that is:
\begin{equation}\label{bias-purf}
 \E[(\tilde{s}(X) - s(X))^2] \leq \frac{6MC^2}{(k+1)^2}
\end{equation}
\begin{flushright}$ \blacksquare $\end{flushright}
\end{prop}

\subsection{Risk bounds for a forest}

Putting together (\ref{var-purf}) and (\ref{bias-purf}) leads to the following risk bound for a PURF.
\begin{theo}\label{theo-risk-forest}
 If $k \xrightarrow[n \to +\infty]{} +\infty$, $\frac{k}{n} \xrightarrow[n \to +\infty]{} 0$, $0 < \mu \leq M$, $s$ is $C$-Lipschitz and $ q \xrightarrow[n \to +\infty]{} +\infty $, the risk of a PURF satisfies:
 $$\E[(\hat{s}(X) - s(X))^2] \leq \frac{3}{4} \frac{\sigma^2 (k+1)}{n} + \frac{6MC^2}{(k+1)^2} + \petito{n \to +\infty} \left( \frac{k}{n} \right)$$
\begin{flushright}$ \blacksquare $\end{flushright}
\end{theo}
Again, taking $(k+1)=n^{1/3}$ gives the upper bound for the risk:
\begin{cor} Under the assumptions of Theorem~\ref{theo-risk-forest},
$$ \E[(\hat{s}(X) - s(X))^2]  \leq K n^{-2/3}  + \petito{n \to +\infty} (n^{-2/3})$$
where $K$ is a positive constant.
\begin{flushright}$ \blacksquare $\end{flushright}
\end{cor}
So, a PURF reaches the minimax rate of convergence for $C$-Lipschitz functions.

Secondly, as the variance of a PUR Forest is systematically reduced compared to a PUR Tree and the bias of a PUR Forest is not larger than the one of a PUR Tree, the risk of a PUR Forest is actually lower.

\section{Conclusion}\label{conclusion}

We emphasize, for a very simple version of random forests, the actual gain of using a random forest instead of using a single random tree. First, we showed that both trees and forests reach the minimax rate of convergence. Then, we manage to highlight a reduction of the variance of a forest, compared to the variance of a tree. This is, in this specific context, a proof of the well-known conjecture for random forests: ``a random forest, by aggregating several random trees, reduces variance and leaves the bias unchanged'' which can be found for example in Hastie et al. \cite{Hastie01}.

An interesting open problem would be to generalize this result, which could handle more complex versions of random forests and relax the hypotheses we made here. Obviously, a more ambitious goal would be to give some precise insights explaining the outstanding performances of random forests-RI.

\section{Proofs}\label{proofs}

\subsection{Proof of Proposition~\ref{prop-var-tree}}\label{proof-prop-var-tree}

We must show that the three last terms in the sum of equation~(\ref{var-tree-gen}) are negligible compared to the constant term $\sigma^2$.

\noindent Let us fix $0 \leq j \leq k$. As it can be found e.g. in Chapter $6$ of \cite{David}, the probability density function of $U_{(j+1)} - U_{(j)}$ is the function $t\in[0,1] \longmapsto k (1-t)^{k-1} $.
\begin{itemize}
 \item For the second term $\E[\delta_{n,p_j}]$:

from \cite{Arlot} we have $\delta_{n,p_j} \leq \kappa_3 (np_j)^{-1/4}$, where $\kappa_3$ is a positive constant. So,
\begin{align*}
\E[ \delta_{n,p_j} ] & \leq  \kappa_3 \E[ (np_j)^{-1/4} ]\\
		     & =     \frac{\kappa_3}{n^{-1/4}} \E[ p_j^{-1/4} ]\\
		     & \leq  \frac{\kappa_3}{(mn)^{-1/4}} \E[ ( U_{(j+1)} - U_{(j)} )^{-1/4} ]\\
		     & \leq  \frac{\kappa_4}{m^{-1/4}} \left(  \frac{k}{n} \right) ^{1/4}
\end{align*}
where $m=\min_{[0,1]} \mu$ and $\kappa_4$ is another positive constant.

Since $\frac{k}{n} \xrightarrow[n \to +\infty]{} 0$ the last upper bound tends to $0$ as $n$ tends to infinity.

 \item For the third term $\E[(\sigma_j^d)^2]$:

$\begin{array}{lll}
(\sigma_j^d)^2 & = & \E[(s(X)- \tilde{s}_{\mathbb{U}}(X))^2 \, | \, U_{(j)} < X \leq U_{(j+1)}]\\
	       & \leq & C^2 ( U_{(j+1)} - U_{(j)} )^2 \quad \mbox{because } s\mbox{ is } C-\mbox{Lipschitz}
\end{array}
$.

So, $\E[\sigma_j^d)^2] \leq C^2 \E[ ( U_{(j+1)} - U_{(j)} )^2 ] = C^2 \frac{2}{(k+1)(k+2)}$ which tends to $0$ as $k$ tends to infinity.

 \item For the last term, the following inequality is sufficient to conclude:

$ \E[(\sigma_j^d)^2\delta_{n,p_j}] \leq C^2\E[\delta_{n,p_j}]$, because $ U_{(j+1)} - U_{(j)} \leq 1 $.

\end{itemize}

\subsection{Proof of Proposition~\ref{prop-bias-tree}}\label{proof-prop-bias-tree}

Function $s$ is supposed to be $C$-Lipschitz, so
\begin{align*}
\E[(\tilde{s}_{\mathbb{U}}(X) - s(X))^2] & =  \E[(\sum_{j=0}^k (s(X)-\beta_j) \mathds{1}_{U_{(j)} < X \leq U_{(j+1)}} )^2] \\
 & =  \E[\sum_{j=0}^k (s(X)-\beta_j)^2 \; \mathds{1}_{U_{(j)} < X \leq U_{(j+1)}} ] \\
 & \leq  \E[\sum_{j=0}^k C^2 ( U_{(j+1)} - U_{(j)} )^2 \; \mathds{1}_{U_{(j)} < X \leq U_{(j+1)}} ] \\
 & =  C^2 \E[\sum_{j=0}^k  ( U_{(j+1)} - U_{(j)} )^2 \; \P(U_{(j)} < X \leq U_{(j+1)}) ] \\
 & \leq  C^2 \E[\sum_{j=0}^k M ( U_{(j+1)} - U_{(j)} )^3] \\
 & \hspace{2cm} \mbox{ because }\mu \mbox{ is bounded by } M\\
 & =  MC^2 \sum_{j=0}^k \E[(U_{(j+1)} - U_{(j)} )^3 ] \\
 & =  MC^2 \frac{6}{(k+2)(k+3)} \\
 & \leq  \frac{6MC^2}{(k+1)^2} \; .
\end{align*}

\subsection{Proof of Theorem~\ref{theo-var-forest}}\label{proof-theo-var-forest}

Before entering into details of the proof of Theorem~\ref{theo-var-forest}, we recall that in the proof of Proposition~\ref{Arlot} (which can be found in \cite{Arlot}), calculations lead to the following equality:
\begin{equation}\label{variance}
 \E[(\hat{s}_{\mathbb{U}}(X) - \tilde{s}_{\mathbb{U}}(X))^2 \, | \, \mathbb{U}] = \sum_{j=0}^k p_j \E \Big[ \frac{1}{n \hat{p}_j} \Big](\sigma^2 + (\sigma_j^d)^2)
\end{equation}
\noindent where $\hat{p}_j=\frac{\sharp\{i: \; U_{(j)} < X_i \leq U_{(j+1)} \}}{n}$.

\noindent Then, an estimation of $p_j \E \Big[ \frac{1}{n \hat{p}_j} \Big]$ gives the expression $\frac{1}{n} (1+\delta_{n,p_j})$ in Proposition~\ref{Arlot}.

\noindent We note
\begin{equation}\label{elementary-variance}
Var_j=p_j \E \Big[ \frac{1}{n \hat{p}_j} \Big](\sigma^2 + (\sigma_j^d)^2)
\end{equation}
a generic term of the sum in the r.h.s. of (\ref{variance}).

\bigskip

We now address the proof of Theorem~\ref{theo-var-forest}. We begin by introducing some notations and establish an intermediate result. The following proposition is not only useful to prove Theorem~\ref{theo-var-forest}, but has its own interest. Indeed, it gives a general upper bound (to be compared to equation~(\ref{var-cond})) which does not depend on the distribution of random partitions defining the trees.

\bigskip

In the sequel we denote the covariance between two PURT by:
\begin{equation*}
 \C(\hat{s}_{\mathbb{U}^1},  \hat{s}_{\mathbb{U}^2}) = \E[(\hat{s}_{\mathbb{U}^1}(X) - \tilde{s}_{\mathbb{U}^1}(X))(\hat{s}_{\mathbb{U}^2}(X) - \tilde{s}_{\mathbb{U}^2}(X))]
\end{equation*}

Let us consider $\mathbb{U}^1=( U_1^1 , \dots , U_k^1 )$ and $\mathbb{U}^2=( U_1^2 , \dots , U_k^2 )$ two sequences of i.i.d. uniform random variables, with respective ordered statistics $( U_{(1)}^1 , \dots , U_{(k)}^1 )$ and $( U_{(1)}^2 , \dots , U_{(k)}^2 )$.

Then we denote by $( V_{(1)} , \dots , V_{(2k)} )$ the ordered statistics of the complete vector $( U_1^1 , \dots , U_k^1 , U_1^2, \ldots , U_k^2 )$, $V_{(0)}=0 $ and $V_{(2k+1)}=1$.

$(\Sigma_t^{d,1,2})^2$ denotes a sum of terms $ \E[(\tilde{s}_{\mathbb{U}^1}(X) - s(X))(\tilde{s}_{\mathbb{U}^2}(X) - s(X)) \, | \, V_{(t')} < X \leq V_{(t'+1)}] $ for several consecutive values of $t'$.

Finally $ \tilde{p}_t $ denotes for some $j\in\{0,\ldots,k\}$ either $ p_j^1 $ or $ p_j^2 $ depending on the relative positions between the $( U_1^1 , \dots , U_k^1 )$ and the $ ( U_1^2 , \dots , U_k^2 )$ in $( V_{(1)} , \dots , V_{(2k)} )$ (see details below).

\begin{prop}\label{prop-var-forest}
The covariance between two PURT satisfies the following upper bound:
\begin{equation}\label{prop-var-forest-gen}
\C(\hat{s}_{\mathbb{U}^1},  \hat{s}_{\mathbb{U}^2}) \leq \frac{1}{n} \E\left[ \sum_{t=0}^{N_{1,2}} (1+\delta_{n,\tilde{p}_t}) (\sigma^2 + (\Sigma_t^{d,1,2})^2) \right]
\end{equation}
where $N_{1,2} = k+1 - \sum_{r=1}^{k-2} \sum_{s=1}^{k-1} \mathds{1}_{U_{(s)}^2 < U_{(r)}^1 < U_{(r+1)}^1 < U_{(r+2)}^1 < U_{(s+1)}^2}$  .
\begin{flushright}$ \blacksquare $\end{flushright}
\end{prop}

\begin{rem}
The gain in variance for a PURF comes from the fact that the number of terms in the sum of equation~(\ref{prop-var-forest-gen}) is smaller than $k+1$. Indeed, it is $k+1-M_{1,2}$ where $M_{1,2}$ is the number of times that $3$ consecutive ordered statistics of $\mathbb{U}^1$ are included in $2$ consecutive ordered statistics of $\mathbb{U}^2$.
\end{rem}

We now prove inequality~(\ref{prop-var-forest-gen}) of Proposition~\ref{prop-var-forest}.
\noindent The term $(\hat{s}_{\mathbb{U}^1}(X) - \tilde{s}_{\mathbb{U}^1}(X))(\hat{s}_{\mathbb{U}^2}(X) - \tilde{s}_{\mathbb{U}^2}(X))$ equals, by definition, to:
\begin{align}\label{covariance}
& \left( \sum_{r=0}^k (\hat{\beta}_r^1 - \beta_r^1) \mathds{1}_{U_{(r)}^1 < X \leq U_{(r+1)}^1} \right) \left( \sum_{s=0}^k (\hat{\beta}_s^2 - \beta_s^2) \mathds{1}_{U_{(s)}^2 < x \leq U_{(s+1)}^2} \right) \notag \\
& =\sum_{t=0}^{2k} (\hat{\beta}_{t,r}^1 - \beta_{t,r}^1) (\hat{\beta}_{t,s}^2 - \beta_{t,s}^2) \mathds{1}_{ V_{(t)} < X \leq V_{(t+1)} }
\end{align}
\noindent where $( V_{(1)} , \dots , V_{(2k)} )$ is the ordered statistics of the vector \\
$( U_1^1 , \dots , U_k^1 , U_1^2, \ldots , U_k^2 )$, $V_{(0)}=0 $, $V_{(2k+1)}=1 $, and

\medskip

$\left\{
\begin{array}{cc}

\hat{\beta}_{t,r}^1 = \hat{\beta}_r^1 \mbox{ and } \beta_{t,r}^1 = \beta_r^1, & \mbox{if } ] V_{(t)} , V_{(t+1)} ] \subset ] U_{(r)}^1 , U_{(r+1)}^1 ] \\
\hat{\beta}_{t,s}^2 = \hat{\beta}_s^2 \mbox{ and } \beta_{t,s}^2 = \beta_s^2, & \mbox{if } ] V_{(t)} , V_{(t+1)} ] \subset ] U_{(s)}^2 , U_{(s+1)}^2 ] \\
\end{array}
\right.$

\medskip

\noindent For $ l=1,2 $ and $ j=0,\ldots,k $, we define $ \hat{p}_j^l = \frac{\sharp\{i: \; U_{(j)}^l < X_i \leq U_{(j+1)}^l \}}{n} $.

Now, let us give some details for the first term of (\ref{covariance}), denoted by $S_1(X)$. Without loss of generality, we suppose that $V_{(1)} = U_{(1)}^1$ (i.e. $U_{(1)}^1 < U_{(1)}^2$). So,
\begin{align*}
 S_1(X) & =  (\hat{s}_{\mathbb{U}^1}(X) - \tilde{s}_{\mathbb{U}^1}(X))(\hat{s}_{\mathbb{U}^2}(X) - \tilde{s}_{\mathbb{U}^2}(X)) \mathds{1}_{0 < X \leq U_{(1)}^1} \\
& =  (\hat{\beta}_1^1 - \beta_1^1) (\hat{\beta}_1^2 - \beta_1^2) \mathds{1}_{0 < X \leq U_{(1)}^1} \\
& =  \Big( \frac{1}{n \hat{p}_1^1} \sum_{i: \; 0 < X_i \leq U_{(1)}^1} (Y_i - \beta_1^1) \Big) \Big( \frac{1}{n \hat{p}_1^2} \sum_{i: \; 0 < X_i \leq U_{(1)}^2} (Y_i - \beta_1^2) \Big) \mathds{1}_{0 < X \leq U_{(1)}^1} \\
& =  \frac{1}{n \hat{p}_1^1 n \hat{p}_1^2} \sum_{i^1: \; 0 < X_{i^1} \leq U_{(1)}^1 \newline i^2: \; 0 < X_{i^2} \leq U_{(1)}^2} (Y_{i^1} - \beta_1^1) (Y_{i^2} - \beta_1^2) \mathds{1}_{0 < X \leq U_{(1)}^1} \\
\end{align*}
\noindent If we denote by $\E^{\Lambda_{1,2}} [.] $ the conditional expectation\\
$\E[.\, | \, (\mathds{1}_{0 < X_{i^1} \leq U_{(1)}^1})_{1 \leq i^1 \leq n} , (\mathds{1}_{0 < X_{i^2} \leq U_{(1)}^2})_{1 \leq i^2 \leq n} ] $, we have:

\medskip

\noindent $ \E[S_1(X) \, | \, \mathbb{U}^1 , \mathbb{U}^2] \newline
= \E \Big[ p_1^1 \E \Big[  \frac{1}{n \hat{p}_1^1 n \hat{p}_1^2} \sum_{i^1: \; 0 < X_{i^1} \leq U_{(1)}^1 \\ i^2: \; 0 < X_{i^2} \leq U_{(1)}^2} \E^{\Lambda_{1,2}}[ (Y_{i^1} - \beta_1^1) (Y_{i^2} - \beta_1^2) ] \Big] \: \Big| \, \mathbb{U}^1 , \mathbb{U}^2 \Big]
$

\medskip

\noindent but
$$i^1 \neq i^2 \quad \Longrightarrow \quad \E^{\Lambda_{1,2}}[ (Y_{i^1} - \beta_1^1) (Y_{i^2} - \beta_1^2) ] =0 $$
because $Y_{i^1}$ and $Y_{i^2}$ are independent. Hence:
\begin{align*}
& \E[S_1(X) \, | \, \mathbb{U}^1 , \mathbb{U}^2] \\
& = \E \Big[ p_1^1 \E \Big[  \frac{1}{n \hat{p}_1^1 n \hat{p}_1^2} \sum_{i: \; 0 < X_i \leq U_{(1)}^1} \E^{\Lambda_1}[ (Y_i - \beta_1^1) (Y_i - \beta_1^2) ]  \Big]   \Big| \, \mathbb{U}^1 , \mathbb{U}^2  \Big] \\
& = \E \Big[ p_1^1 \E \Big[  \frac{1}{n \hat{p}_1^1 n \hat{p}_1^2} \sum_{i: \; 0 < X_i \leq U_{(1)}^1} \E[ (Y_i - \beta_1^1) (Y_i - \beta_1^2) \, | \, 0 < X_i \leq U_{(1)}^1 ]  \Big]   \Big| \, \mathbb{U}^1 , \mathbb{U}^2  \Big]
\end{align*}

\noindent where $\E^{\Lambda_1} [.] $ denotes the conditional expectation $\E[.\, | \, (\mathds{1}_{0 < X_i \leq U_{(1)}^1})_{1 \leq i \leq n} ]$.

\bigskip

\noindent Now, as
$$\E[ (Y_i - \beta_1^1) (Y_i - \beta_1^2) \, | \, 0 < X_i \leq U_{(1)}^1 ] = \E[ (Y - \beta_1^1) (Y - \beta_1^2) \, | \, 0 < X \leq U_{(1)}^1 ]$$
for all $i$, and
$$ \E[ (Y - \beta_1^1) (Y - \beta_1^2) \, | \, 0 < X \leq U_{(1)}^1 ] = \sigma^2 + ( \sigma_0^{d,1,2} )^2 $$
where
$$ ( \sigma_0^{d,1,2} )^2 = \E[ (s(X) - \tilde{s}_{\mathbb{U}^1}(X)) (s(X) - \tilde{s}_{\mathbb{U}^2}(X)) \, | \, 0 < X \leq V_{(1)} ] $$
we get
$$ \E[S_1(X) \, | \, \mathbb{U}^1 , \mathbb{U}^2] = p_1^1 \E \Big[ \frac{1}{n \hat{p}_1^2}  \Big] ( \sigma^2 + ( \sigma_0^{d,1,2} )^2 ) \; . $$

\medskip

If we suppose in addition that $V_{(2)} = U_{(1)}^2 $, we similarly get for the second term of (\ref{covariance}):

\medskip

\noindent $ \E[S_2(X) \, | \, \mathbb{U}^1 , \mathbb{U}^2] \newline
=\E[(\hat{s}_{\mathbb{U}^1}(X) - \tilde{s}_{\mathbb{U}^1}(X))(\hat{s}_{\mathbb{U}^2}(X) - \tilde{s}_{\mathbb{U}^2}(X)) \mathds{1}_{U_{(1)}^1 < X \leq U_{(1)}^2} \, | \, \mathbb{U}^1 , \mathbb{U}^2] \newline
= q_2 \E \left[ \frac{n \hat{q}_2}{n \hat{p}_2^1 n \hat{p}_1^2} ( \sigma^2 + ( \sigma_1^{d,1,2} )^2 ) \right] $

\noindent where

\medskip

$$q_2 = \P( V_{(1)} < X \leq V_{(2)}) = \P( U_{(1)}^1 < X \leq U_{(1)}^2 ) $$
$$n \hat{q}_2 = \sharp \{ i: V_{(1)} < X_i \leq V_{(2)} \} $$
and
$$( \sigma_1^{d,1,2} )^2 = \E[ (s(X) - \tilde{s}_{\mathbb{U}^1}(X)) (s(X) - \tilde{s}_{\mathbb{U}^2}(X)) \, | \, V_{(1)} < X \leq V_{(2)} ] \; . $$

\noindent Since $]V_{(1)}, V_{(2)}]$ is included in $]U_{(1)}^1, U_{(2)}^1]$, we have $ \hat{q}_2 \leq \hat{p}_2^1 $, so:

$$ \E[S_2(X) \, | \, \mathbb{U}^1 , \mathbb{U}^2] \leq q_2 \E \Big[ \frac{1}{n \hat{p}_1^2}  \Big] ( \sigma^2 + ( \sigma_1^{d,1,2} )^2 ) \; . $$

\medskip

\noindent Finally, by summing the two terms $S_1(X)$ and $S_2(X)$, we deduce that

\medskip

$\E[S_1(X) + S_2(X) \, | \, \mathbb{U}^1 , \mathbb{U}^2] \leq p_1^2 \E \Big[ \frac{1}{n \hat{p}_1^2}  \Big] ( \sigma^2 + ( \sigma_0^{d,1,2} )^2 + ( \sigma_1^{d,1,2} )^2 ) $

\medskip

\noindent In conclusion, we succeeded to bound the sum of the first two terms of (\ref{covariance}) by an expression very close to $ Var_j $ (defined in (\ref{elementary-variance})). The only difference comes from the fact that instead of $(\sigma_j^d )^2$ we have $( \sigma_0^{d,1,2} )^2 + ( \sigma_1^{d,1,2} )^2$. But as we saw in proof of Proposition~\ref{prop-var-tree}, these terms are negligible, so $p_1^2 \E \Big[ \frac{1}{n \hat{p}_1^2}  \Big] ( \sigma^2 + ( \sigma_0^{d,1,2} )^2 + ( \sigma_1^{d,1,2} )^2 )$ is of the same order than $ Var_j $.

We can easily generalize this fact by proving the following lemma.

We denote by $S_j(X)$ the $j$-th term of (\ref{covariance}), i.e. $S_j(X)=(\hat{s}_{\mathbb{U}^1}(X) - \tilde{s}_{\mathbb{U}^1}(X))(\hat{s}_{\mathbb{U}^2}(X) - \tilde{s}_{\mathbb{U}^2}(X)) \mathds{1}_{V_{(j)} < X \leq V_{(j+1)}}$.

\begin{lem}\label{lemun}
Let $r$ be in $\{0,\ldots,k\}$ and denote by $t$, $t'$ the integers such that
\begin{equation}\label{hyp-lemun}
 U^1_{(r)} = V_{(t)} < V_{(t'+1)} = U^1_{(r+1)}
\end{equation}
then
$$\E \bigg[ \sum_{j=t}^{t'} S_j(X) \, | \, \mathbb{U}^1 , \mathbb{U}^2 \bigg] \leq p_r^1 \E \Big[ \frac{1}{n \hat{p}_r^1}  \Big] ( \sigma^2 + (\Sigma_r^{d,1,2})^2 )$$
where $(\Sigma_r^{d,1,2})^2 = \sum_{j=t}^{t'} ( \sigma_j^{d,1,2} )^2 $.
\begin{flushright}$ \blacksquare $\end{flushright}
\end{lem}

\noindent Indeed for all $j \in \{t,t+1,\ldots,t'\}$,
$$\E[S_j(X) \, | \, \mathbb{U}^1 , \mathbb{U}^2] \leq q_j \E \Big[ \frac{1}{n \hat{p}_r^1}  \Big] ( \sigma^2 + ( \sigma_j^{d,1,2} )^2 )$$
where
$$q_j=\P(V_{(j)} < X \leq V_{(j+1)})$$
and
$$( \sigma_j^{d,1,2} )^2 = \E[ (s(X) - \tilde{s}_{\mathbb{U}^1}(X)) (s(X) - \tilde{s}_{\mathbb{U}^2}(X)) \, | \, V_{(j)} < X \leq V_{(j+1)} ] \; .$$

\noindent Thus,
$$\E \bigg[ \sum_{j=t}^{t'} S_j(X) \, | \, \mathbb{U}^1 , \mathbb{U}^2 \bigg] \leq \P(V_{(t)} < X \leq V_{(t'+1)}) \; \E \Big[ \frac{1}{n \hat{p}_r^1}  \Big] ( \sigma^2 + (\Sigma_r^{d,1,2})^2 ) \; .$$
From relation~(\ref{hyp-lemun}) we have $ \P(V_{(t)} < X \leq V_{(t'+1)}) = p_r^1 $, which concludes the proof of Lemma~\ref{lemun}.

\bigskip

Therefore, we can upper bound the initial sum (\ref{covariance}) of $2k+1$ terms by a sum of $k+1$ terms of the same order as $ Var_j $ only involving intervals of the partition $ \mathbb{U}^1 $. At this stage, we get an upper bound for the variance of a forest which is of the same order as the variance of a tree. But we can do better. With similar arguments, we can prove the following lemma:

\begin{lem}\label{lemdeux}
If there exist $r$ and $s$ such as
$$ U_{(s)}^2 < U_{(r)}^1 < U_{(r+1)}^1 < U_{(r+2)}^1 < U_{(s+1)}^2 $$
the expression
$$ \E[(\hat{s}_{\mathbb{U}^1}(X) - \tilde{s}_{\mathbb{U}^1}(X))(\hat{s}_{\mathbb{U}^2}(X) - \tilde{s}_{\mathbb{U}^2}(X)) \mathds{1}_{U_{(r)}^1 < X \leq U_{(r+2)}^1} \, | \, \mathbb{U}^1 , \mathbb{U}^2]  $$
is upper bounded by
$$ p_s^2 \E \Big[ \frac{1}{n \hat{p}_s^2} \Big] ( \sigma^2 + ( \Sigma_{s}^{d,1,2} )^2 ) \; . $$
where $( \Sigma_{s}^{d,1,2} )^2 = ( \sigma_{r+s}^{d,1,2} )^2 + ( \sigma_{r+s+1}^{d,1,2} )^2 $.
\begin{flushright}$ \blacksquare $\end{flushright}
\end{lem}

\noindent Indeed,
\begin{align*}
& \E[(\hat{s}_{\mathbb{U}^1}(X) - \tilde{s}_{\mathbb{U}^1}(X))(\hat{s}_{\mathbb{U}^2}(X) - \tilde{s}_{\mathbb{U}^2}(X)) \mathds{1}_{U_{(r)}^1 < X \leq U_{(r+1)}^1} \, | \, \mathbb{U}^1 , \mathbb{U}^2] \\
& \leq p_r^1 \E \Big[ \frac{1}{n \hat{p}_s^2}  \Big] ( \sigma^2 + ( \sigma_{r+s}^{d,1,2} )^2 )
\end{align*}
and
\begin{align*}
& \E[(\hat{s}_{\mathbb{U}^1}(X) - \tilde{s}_{\mathbb{U}^1}(X))(\hat{s}_{\mathbb{U}^2}(X) - \tilde{s}_{\mathbb{U}^2}(X)) \mathds{1}_{U_{(r+1)}^1 < X \leq U_{(r+2)}^1} \, | \, \mathbb{U}^1 , \mathbb{U}^2]  \\
& \leq p_{r+1}^1 \E \Big[ \frac{1}{n \hat{p}_s^2}  \Big] ( \sigma^2 + ( \sigma_{r+s+1}^{d,1,2} )^2 ) \; .
\end{align*}
Finally, since $p_r^1+p_{r+1}^1 \leq p_s^2$, $( \sigma_{r+s}^{d,1,2} )^2 \leq ( \sigma_{r+s}^{d,1,2} )^2 + ( \sigma_{r+s+1}^{d,1,2} )^2 $ and $( \sigma_{r+s+1}^{d,1,2} )^2 \leq ( \sigma_{r+s}^{d,1,2} )^2 + ( \sigma_{r+s+1}^{d,1,2} )^2 $, the result is obtained by summing the two terms.

\bigskip

As in Proposition~\ref{Arlot}, we replace all $ p_{j}^l \E \Big[ \frac{1}{n \hat{p}_j^l}  \Big] $ by their estimates $ (1+\delta_{n p_j^l})$.

\bigskip

By repeatedly applying this lemma for all intervals, we can upper bound
$$ \E[ (\hat{s}_{\mathbb{U}^1}(X) - \tilde{s}_{\mathbb{U}^1}(X))(\hat{s}_{\mathbb{U}^2}(X) - \tilde{s}_{\mathbb{U}^2}(X)) \, | \, \mathbb{U}^1 , \mathbb{U}^2]$$
by a sum of $N_{1,2}$ terms of the form $(1+\delta_{n,\tilde{p}_t}) (\sigma^2 + (\Sigma_t^{d,1,2})^2)$, where $ \tilde{p}_t $ denotes for some $j\in\{0,\ldots,k\}$ either $ p_j^1 $ or $ p_j^2 $ depending on the fact that we are in the situation of Lemma~\ref{lemun} or Lemma~\ref{lemdeux}, $ N_{1,2} = k+1 - M_{1,2} $ and
$$ M_{1,2} = \sum_{r=1}^{k-2} \sum_{s=1}^{k-1} \mathds{1}_{U_{(s)}^2 < U_{(r)}^1 < U_{(r+1)}^1 < U_{(r+2)}^1 < U_{(s+1)}^2} \; .$$
This concludes the proof of Proposition~\ref{prop-var-forest}. Now, using the fact that we deal with uniform partitions, we manage to prove the following corollary.

\begin{cor}\label{coro}
If $k \xrightarrow[n \to +\infty]{} +\infty$, $\frac{k}{n} \xrightarrow[n \to +\infty]{} 0$, $\mu>0$ and \\
$s$ is $C$-Lipschitz,  we have,
\begin{align*}
\C(\hat{s}_{\mathbb{U}^1},  \hat{s}_{\mathbb{U}^2}) & \leq  \frac{\sigma^2 \E[N_{1,2}]}{n} + \petito{n\to+\infty} \left( \frac{k}{n} \right) \\
 & \leq  \frac{3}{4}\frac{\sigma^2(k+1)}{n} + \petito{n\to+\infty} \left( \frac{k}{n} \right) \; .
\end{align*}
\begin{flushright}$ \blacksquare $\end{flushright}
\end{cor}

Because of the simple draws of random partitions, the number $M_{1,2}$ is explicitly computable (we know the distribution of the two ordered statistics) and it is shown to be equivalent to $\frac{1}{4} (k+1)$ as $k$ tends to $+\infty$ (see Lemma~\ref{sum} below).

\noindent As in Proposition~\ref{prop-var-tree}, we have to prove that all terms of the sum are negligible compared to the constant one $\sigma^2$. To deal with the fact that the number of terms in the sum is now random, we use the following simple inequality:

\medskip

\noindent $ \E\left[ \sum_{t=0}^{N_{1,2}} ( \sigma^2\delta_{n,p_t} + (\Sigma_t^{d,1,2})^2 + (\Sigma_t^{d,1,2})^2 \delta_{n,p_t}) \right]\\
\leq \sum_{t=0}^k \left( \E[\sigma^2\delta_{n,p_t}] + \E[(\Sigma_t^{d,1,2})^2] + \E[(\Sigma_t^{d,1,2})^2 \delta_{n,p_t}] \right).$

\medskip

\noindent These quantities are of the same kind as the three last terms in the sum of equation~\ref{var-tree-gen}. So with the same techniques we get that
$$\frac{1}{n} \E\left[ \sum_{t=0}^{N_{1,2}} ( \sigma^2\delta_{n,p_t} + (\Sigma_t^{d,1,2})^2 + (\Sigma_t^{d,1,2})^2 \delta_{n,p_t}) \right] = \petito{n\to +\infty} \left( \frac{k}{n} \right) \; . $$

\noindent So, we have
$$\E[(\hat{s}_{\mathbb{U}^1}(X) - \tilde{s}_{\mathbb{U}^1}(X))(\hat{s}_{\mathbb{U}^2}(X) - \tilde{s}_{\mathbb{U}^2}(X))]
  \leq  \frac{\sigma^2 \E[N_{1,2}]}{n} + \petito{n\to+\infty} \left( \frac{k}{n} \right)  \; . $$

\noindent Finally, the following technical result allows to conclude the proof of Corollary~\ref{coro}, and thus the proof of Theorem\ref{theo-var-forest}.

\begin{lem}\label{sum}
 $$\E[M_{1,2}] = \frac{(k-2)(k-3)}{2(2k-1)} \left( 1+\frac{4}{(k+1)(k-3)} \right) \; .$$
\noindent Hence,
$$\E[M_{1,2}] = \frac{k+1}{4} + \petito{k\to+\infty} (k) \; .$$
\begin{flushright}$ \blacksquare $\end{flushright}
\end{lem}

\noindent We then obtain that
$$ \E[N_{1,2}] = \frac{3}{4} (k+1) + \petito{k\to+\infty} (k)  \; .$$

\noindent Let us demonstrate lemma~\ref{sum}.

$$ \E[M_{1,2}] = \sum_{r=1}^{k-2} \sum_{s=1}^{k-1} \P(U_{(s)}^2 < U_{(r)}^1 < U_{(r+1)}^1 < U_{(r+2)}^1 < U_{(s+1)}^2) $$

\noindent As we know the distribution of ordered statistics (see e.g. Section $2.2$ of \cite{David}), we can compute the following probability:
\begin{align*}
&\P(U_{(s)}^2 < U_{(r)}^1 < U_{(r+1)}^1 < U_{(r+2)}^1 < U_{(s+1)}^2) \\
&= \P(U_{(s)}^2 < U_{(r)}^1 \mbox{ and } U_{(r+2)}^1 < U_{(s+1)}^2) \\
&= \sum_{j=r+2}^k \sum_{i=0}^{r-1} \frac{k!}{i! (j-i)! (k-j)!} \E[(U^2_{(s)})^i (U^2_{(s+1)} - U^2_{(s)})^{j-i} (1 - U^2_{(s+1)})^{k-j}] \\ 
&= \sum_{j=r+2}^k \sum_{i=0}^{r-1} \frac{k!}{i! (k-j)!} \frac{k!}{(s-1)! (k-(s+1))!}\frac{(i+s-1)! (2k-(j+s)-1)!}{(2k)!} 
\end{align*}
\noindent So,
\begin{align*}
& \E[M_{1,2}] \\
& = \frac{(k!)^2}{(2k)!} \sum_{r=1}^{k-2} \sum_{s=1}^{k-1} \left( \sum_{i=0}^{r-1} \left( \begin{array}{c} i+(s-1)\\ i \end{array} \right) \right)
\left( \sum_{j=r+2}^k \left( \begin{array}{c} k-j+k-(s+1)\\ k-j \end{array} \right) \right) \\
& = \frac{(k!)^2}{(2k)!} \sum_{r=1}^{k-2} \sum_{s=1}^{k-1} \left( \begin{array}{c} r-1+s\\ r-1 \end{array} \right)
\left( \begin{array}{c} 2k-r-2-s\\ k-r-2 \end{array} \right)
\end{align*}
(by elementary properties of binomial coefficients (see e.g. \cite{Graham} p.$160$))

\medskip

\noindent $= \frac{k-2}{4(2k-1)} \sum_{t=0}^{2k-5} \sum_{r=t-k+2}^t
\frac{\left( \begin{array}{c} t+1\\ r \end{array} \right)
\left( \begin{array}{c} 2k-3-(t+1)\\ k-3-r \end{array} \right) }
{\left( \begin{array}{c} 2k-3\\ k-3 \end{array} \right)}$

\medskip

(by defining $t=r+s$)

\medskip

\noindent $ = \frac{k-2}{4(2k-1)} \sum_{t=0}^{2k-5}[\mathrm{F}_{\mathcal{H}(2k-3,t+1,k-3)} (t) - \mathrm{F}_{\mathcal{H}(2k-3,t+1,k-3)} (t-k+1)]$

\medskip

(where $\mathrm{F}_{\mathcal{H}(N,m,n)}$ denotes the cumulative distribution function of the hyper-geometric distribution)

\medskip

\noindent $ = \frac{k-2}{4(2k-1)} 2 \sum_{t=0}^{k-3} \mathrm{F}_{\mathcal{H}(2k-3,t+1,k-3)} (t) \newline
= \frac{k-2}{2(2k-1)} \left[ \sum_{t=0}^{k-4} \left( 1 -
\frac{\left( \begin{array}{c} t+1\\ t+1 \end{array} \right)
\left( \begin{array}{c} 2k-3-(t+1)\\ k-3-(t+1) \end{array} \right) }
{\left( \begin{array}{c} 2k-3\\ k-3 \end{array} \right)}
 \right) +1 \right] \newline
= \frac{k-2}{2(2k-1)} \Big( k-3 + \frac{4}{k+1} \Big)
$.

\end{document}